\documentclass{amsart}
\usepackage{amsfonts,latexsym,amsthm,amssymb,xspace}

\newcommand{\mes}{\mathop{\mathrm{mes}}}
\newcommand{\N}{\mathbb{N}}

\newtheorem{thm}{Theorem}
\newtheorem{lem}[thm]{Lemma}
\newtheorem{rem}[thm]{Remark}
\newtheorem{cor}[thm]{Corollary}
\newtheorem{defi}{Definition}

\newcommand{\R}[1][n]{\ensuremath{\mathbb{R}^{#1}}}

\begin{document}
\sloppy

\title{Kruglov operator and operators defined by random permutations}

\author{S.V.~Astashkin${}^1$}
\thanks{${}^1$Research is partially supported by Russian Foundation for Basic Research, 07-01-96603.}
\author{D.V.~Zanin}
\author{E.M.~Semenov${}^2$}
\thanks{${}^2$Research is partially supported by Russian Foundation for Basic Research, 05-01-00629.}
\author{F.A.~Sukochev${}^3$}
\thanks{${}^3$Research is partially supported by Australian Research Council.}

\begin{abstract}
The Kruglov property and the Kruglov operator play an important role
in the study of geometric properties of r.i. function spaces. We prove that
the boundedness of the Kruglov operator in a r.i. space is
equivalent to the uniform boundedness on this space of a sequence of operators defined by
random permutations. It is shown also that there is no minimal r.i.
space with the Kruglov property.
\end{abstract}

\maketitle


\section{Introduction}

Let $f$ be a random variable (measurable function) on the interval $[0,1]$. We denote a random variable $\sum_{i=1}^N f_i$ by $\pi(f).$ Here, $f_i$'s are independent copies of $f$ and $N$ is a Poisson random variable with parameter 1, independent from $f_i$'s.

\begin{defi} A r.i. function space $E$ on the interval $[0,1]$ is said to have the Kruglov property ($E\in\mathbb{K}$) if  $f\in E\Longleftrightarrow \pi(f)\in E.$
\end{defi}

This property was introduced and studied by Braverman \cite{lit-2}, exploiting some constructions and ideas from the article \cite{lit-1} by Kruglov. An operator approach to the study of this property was introduced in \cite{lit-3} (see also \cite{lit-12}).

Let $\{B_n\}_{n=1}^\infty$ be a sequnce of mutually disjoint measurable subsets of $[0,1]$ and let $\mes B_n=\dfrac1{en!}$. If $f\in L_1[0,1],$ then set
$$Kf(\omega_0,\omega_1,\ldots)=\sum_{n=1}^\infty \sum_{k=1}^nf(\omega_k)\chi_{B_n}(\omega_0).$$
Here and everywhere else we denote the characteristic function of the set $B$ by $\chi_B.$ It then follows that $K:\,L_1[0,1]\to L_1(\Omega,P)$ is a positive linear operator. Here $(\Omega,P)=\prod_{n=0}^\infty([0,1],\mes)$, where $\mes$ is the Lebesgue measure on $[0,1]$. Since $Kf$ is equidistributed with $\pi(f)$ (see \cite{lit-3}), we may consider $Kf$ as an explicit representation of $\pi(f).$ In particular, an r.i. space $E\in\mathbb{K}$ if and only if $K$ (boundedly) maps $E$ into $E(\Omega,P)$ (see \cite{lit-3}).

We will also use an equivalent representation of the operator $K$ introduced in \cite{lit-3}. Let $f^*$ be decreasing rearrangement of $|f|,$ that is,  $f^*(t)$ decreases on $[0,1]$ and is equimeasurable with $|f(t)|.$ If $f\in L_1[0,1]$ and if $\{B_n\}$ is the same sequence of subsets of $[0,1]$ as above, then let $f_{n,1},f_{n,2},\ldots,f_{n,n}$, $\chi_{B_n}$ be the set of independent functions for every $n\in\N,$ such that $f_{n,k}^*=f^*$ for every $n\in\N$ and $k=1,2,\ldots,n.$ Under these conditions, $Kf(t)$ is defined as a rearrangement of the function
\begin{equation}\label{40-(-1)}
 \sum_{n=1}^\infty \sum_{k=1}^n f_{n,k}(t)\chi_{B_n}(t)\;\;(0\le t\le 1).
\end{equation}
It follows from the definition of an r.i. space and that of the operator $K$ that $\|Kf\|_E\geq e^{-1}\|f\|_E$ for every r.i. space $E$ and for every $f\in E$ (see also \cite[1.6,p.11]{lit-2}). It is shown in \cite{lit-3} that the operator $K$ plays an important role in estimating the norm of sums of  independent random variables through the norm of sums of their disjoint copies. In particular, in \cite{lit-3} the well-known results of Johnson and Schechtman from \cite{lit-4} have been strengthened.

It is well known \cite{lit-1},\cite{lit-2} that the Orlicz space $\exp L_1$ defined by the function $e^t-1$ satisfies the Kruglov property. The latter property also holds for its separable part $(\exp L_1)_0$. Indeed, since $K$ is bounded in $\exp L_1,$ we have  $K((\exp L_1)_0)\subset\overline{K(L_{\infty})}$ (the closure is taken with respect to the norm in $\exp L_1$). However, $K(L_{\infty})\subset(\exp L_1)_0$  \cite[Theorem 4.4]{lit-3}. Since $(\exp L_1)_0$ is a closed subset of $\exp L_1,$ we conclude that the operator $K$ maps $(\exp L_1)_0$ into itself. All previously known r.i. spaces $E$ with the Kruglov property satisfied the inclusion $E\supset (\exp L_1)_0.$ This together with some results from \cite{lit-3} (e.g. Theorem 7.2) suggest that $(\exp L_1)_0$ is the minimal r.i. space with the Kruglov property. However, in the first part of the paper we show that this conjecture fails. Moreover, we show that for every given r.i. space $E\in\mathbb{K}$ there exists a Marcinkiewicz space satisfying the Kruglov property such that $M_\psi\subsetneq E$ (see Corollary \ref{corol: there isn't minim. r.i.s.}). The situation is quite different in the subclass of Lorentz spaces. Indeed, every Lorentz space satisfying the Kruglov property necessarily contains $\exp L_1$ (see Theorem \ref{Lorentz theorem}).

In \cite{lit-5}, Kwapien and Sch\" utt considered random permutations and applied their results to the geometry of Banach spaces. These results were further strengthened in \cite{lit-6} and \cite{lit-7} via an operator approach. The following family of operators was introduced there. Let $n\in\mathbb{N}$ and let $S_n$ be the set of all permutations of scalars $1,2,\cdots,n.$ From now on the sets $S_n$ and $\{1,2,\cdots,n!\}$ will be identified (in an arbitrary manner). Firstly, we define an operator $A_n$ acting from $\mathbb{R}^n$ into $\mathbb{R}^{n!}$: if $x=(x_1,x_2,\dots,x_n)\in\mathbb{R}^n$ and if $\pi\in S_n$ is an arbitrary permutation, then
\begin{equation}\label{40-1}
A_nx(\pi):=\sum_{i:\,\pi(i)=i}x_i.
\end{equation}
For every $x\in L_1[0,1]$, we define a vector $B_nx\in\mathbb{R}^n$ with coordinates $(B_nx)_i=n\int_{(i-1)/n}^{i/n} x(t)\,dt,$ $i=1,2,\dots,n.$ The operator $B_n$ has a right inverse operator $C_n$ ( $B_nC_nx=x$ for every $x\in \mathbb{R}^n$)  which maps every vector into a function with constancy intervals $[(i-1)/n,i/n].$ Now, we define
$$T_n=C_{n!}A_nB_n.$$
For every $n\in\mathbb{N},$ $T_n$ is a positive linear operator from $L_1[0,1]$ into the space of step functions. It is not hard to show that
\begin{equation}\label{40-0}
 \|T_nx\|_{L_1}=\|x\|_{L_1}
\end{equation}
for every positive $x\in L_1[0,1].$ Sometimes, we will also use the notation $T_n$ for the operator $C_{n!}A_n,$ defined analogously on $\mathbb{R}^n$ (this does not cause any ambiguity). If $x=(x_1,x_2,\ldots,x_n)\in\mathbb{R}^n$ and if $E$ is an r.i. space, then the notation $\|x\|_E$ will always mean
$$\|C_nx\|_E=\displaystyle\Big\| \sum_{k=1}^n x_k\chi_{\Big(\frac{k-1}
n,\frac kn\Big)} \Big\|_E.$$

The operators generated by random permutations and defined on the set of square matrices were considered in \cite{lit-7}, where it was established that such operators are uniformly bounded if the family of operators $\{T_n\}_{n\ge 1}$ is uniformly bounded. There is no any visible connection between the operators $K$ and $T_n, \ n\ge 1$. Nevertheless, the following interesting fact follows from the comparison of
results in \cite{lit-7} and \cite{lit-3}: the criterion for the boundedness of the operator $K$ in any Lorentz space $\Lambda_{\varphi}$ and that for the uniform boundedness of the family of operators $\{T_n\}_{n\ge 1}$ in $\Lambda_{\varphi}$ coincide.  More precisely, both criteria are equivalent to the following condition
\begin{equation}\label{40-2}
M:= \sup_{0<t\leq1}\frac1{\varphi(t)} \sum_{k=1}^\infty\varphi\left(\frac{t^k}{k!}\right)<\infty.
\end{equation}
It is now natural to ask  whether the boundedness of the operator $K$ in an arbitrary r.i. space $E$ is equivalent to the uniform boundedness of the family of operators $\{T_n\}_{n\ge 1}$ in $E$.  In the second part of this paper we establish that it is indeed the case. The proof is based on combinatorial arguments and is connected with obtaining estimates of corresponding distribution functions. The established equivalence implies some new corollaries for the operator $K$ and operators $T_n, \ n\ge 1$. In particular, Corollary ~\ref{corol-40-8} strengthens Theorem~19 from \cite{lit-7} by showing that the uniform boundedness of the family of operators $\{T_n\}_{n\ge 1}$ in Orlicz spaces $\exp L_p$ is equivalent to the condition $p\le 1$.

The authors thank the referee for comments and suggestions which allowed to simplify the definition of the operator $T_n$, $n\ge 1$ and the proof of Lemma 7 and in general were helpful in improving the final text of this paper.
\vskip 0.2cm

\section{Definitions and notation}

A Banach space $E$ consisting of functions measurable on $[0,1]$ is said to be rearrangement invariant or symmetric (r.i.) if  the following conditions hold
\begin{enumerate}
\item If $|x(t)|\leq|y(t)|$ for a.e. $t\in[0,1]$ and $y\in E,$ then $x\in E$ and $\|x\|_E\leq\|y\|_E.$
\item If functions $x$ and $y\in E$ are equimeasurable, that is
$$\mes\{t\in [0,1]:\,|x(t)|>\tau\}=\mes\{t\in[0,1]:\,|y(t)|>\tau\}\quad (\tau>0),$$
then $x\in E$ and $\|x\|_E=\|y\|_E$.
\end{enumerate}

If $E$ is an r.i. space, then $L_\infty\subset E\subset L_1$ and these inclusions are continuous. Moreover, if $\|\chi_{(0,1)}\|_E=1,$ then $\|x\|_{L_1}\leq\|x\|_E\leq\|x\|_{L_\infty}$ for every $x\in L_\infty.$ For every $\tau>0,$ the dilation operator $\sigma_{\tau}$ defined by $\sigma_\tau x(t):=x(t/\tau)\chi_{[0,1]}(t/\tau)$ $(0\le t\le 1)$ boundedly maps $E$ into itself and $\|\sigma_\tau\|_E\le \max(1,\tau).$

The K\" othe dual space $E^\prime$ consists of all functions $x$ for which the norm
$$ \|x\|_{E^\prime}=\sup_{\|y\|_E\leq1} \int_0^1x(t)y(t)dt$$
is finite. Clearly, $E^{\prime}$ is also an r.i. space. Following \cite[2.a.1]{lit-8}, we assume that either r.i. space $E$ is separable or $E$ coincides with its second K\" othe dual space $E^{\prime\prime}$. In any case, the space $E$ is contained in $E^{\prime\prime}$ as a closed subspace and the inclusion $E\subset E^{\prime\prime}$ is an isometry. If $E$ is separable, then $E^{\prime}$ coincides with its dual space $E^*$. The closure $E_0$ of $L_{\infty}$ in $E$ is called the separable part of $E$. The space $E_0$ is separable provided that $E\ne L_\infty.$

Recall that the weak convergence of distributions of measurable on  $[0,1]$ functions $x_n$ to the distribution of the function $x$ ($x_n\Rightarrow x$) means that for every continuous and bounded on $(-\infty, \infty)$ function $y$ we have$$\lim_{n\to\infty}\int_{-\infty}^\infty y(t)\,d\mes\{s:\,x_n(s)<t\}=
\int_{-\infty}^\infty y(t)\,d\mes\{s:\,x(s)<t\}.$$
If $E$ is an r.i. space, $x_n\in E$ $(n\in\N),$ $\limsup_{n\to\infty}\|x_n\|_E=C<\infty$ and $x_n\Rightarrow x,$ then $x\in E''$ and $\|x\|_{E''}\le C$ \cite[Proposition~1.5]{lit-2}.

The following submajorization defined on $L_1$ plays an important role in the theory of r.i. spaces. We denote $x\prec y$ if
\[
 \int_0^\tau x^*(t)dt\leq \int_0^\tau y^*(t)dt
\]
for all $\tau\in[0,1]$. If $x\prec y$ and $y\in E,$ then $x\in E$ and $\|x\|_E\leq\|y\|_E$. Here and below, $x^*(t)$ is the non-increasing left continuous rearrangement of the function $|x(t)|,$ i.e.
$$
x^*(t)=\inf\{\tau\ge 0:\,\mes\{s\in
[0,1]:\,|x(s)|>\tau\}<t\}\;\;(0<t\le 1).$$

We list below the most important examples of r.i. spaces. Let $M$ be an increasing convex function on $[0,\infty)$ such that $M(0)=0$. By $L_M$ we denote the Orlicz space $L_M$ with the norm
\[
 \|x\|_{L_M}=\inf\left\{ \lambda>0:\,
 \int_0^1M\left(\frac{|x(t)|}\lambda\right)\,dt\leq1\right\}.
\]
Function $M_p(u)=e^{u^p}-1$ is convex if $p\geq1$ and is equivalent to some convex function if $0<p<1$. We denote $L_{M_p}$ by $\exp L_p$.

Let $\varphi(t)$ be an increasing concave function on $[0,1]$ such that $\varphi(0)=0$ and let $\Lambda_\varphi$ be the Lorentz space equipped with a norm
\[
 \|x\|_{\Lambda_\varphi}=\int_0^1x^*(t)\,d\varphi(t).
\]
Similarly, $M_\varphi$ is the Marcinkiewicz space equipped with the norm
\[
 \|x\|_{M_\varphi}=\sup_{0<t\le
 1}\frac{1}{\varphi(t)}{\int_0^tx^*(s)\,ds}.
\]
All facts listed above from the theory of r.i. spaces and more detailed information about this theory may be found in  the books \cite{lit-8}, \cite{lit-9}.

In what follows, ${\rm supp}\,f$ is the support of the function $f$, i.e. the set $\{t:\,f(t)\ne 0\}$. We write $F\asymp G,$ if $C^{-1}F\le G\le CF,$ where $C>0$ is a constant. Finally, $|A|$  denotes the number of elements of a finite set $A$.
 \vskip 0.2cm

{\centering\section{Lorentz and Marcinkiewicz spaces ``near'' $\exp L_1$}}

\begin{thm}\label{psi theorem} There exists a family of Marcinkiewicz spaces $\{M_{\psi_{\varepsilon}}\}_{0<\varepsilon<1}$ such that $M_{\psi_{\varepsilon}}\subset M_{\psi_{\delta}}$ for every $0<\varepsilon\leq\delta<1,$ satisfying the following conditions:
\begin{enumerate}
\item $M_{\psi_{\varepsilon}}\in\mathbb{K}$, $0<\varepsilon<1$.
\item For every r.i. space $E\in\mathbb{K}$ we have  $M_{\psi_{\varepsilon}}\subset E$ if $\varepsilon$ is small enough.
\item Functions $\psi_{\varepsilon}$ are not pairwise equivalent, or more precisely,
\begin{equation}\label{non-equivalent}
\lim_{t\to0}\frac{\psi_{\varepsilon}(t)}{\psi_{\delta}(t)}=0,\quad if\quad 0<\varepsilon<\delta<1.
\end{equation}
\item We have $M_{\psi_{\varepsilon}}\varsubsetneqq(\exp L_1)_0$ if $\varepsilon>0$ is small enough.
\end{enumerate}
\end{thm}

We will need the following simple assertion.

\begin{lem}\label{power of operator K} For every $f\in L_1[0,1]$
$$\lim_{n\to\infty}\mes({\rm supp}\,K^nf)=0.$$
\end{lem}
\begin{proof}
Since the operator $K$ is positive, we may assume that $f\ge 0$ and that $\mes({\rm supp}\,f)=1$.
If $a_n:=\mes\{t:\,K^nf(t)=0\}$ $(n\in\N),$ then, by definition of the operator $K$ (see equation \eqref{40-(-1)})
$a_1=1/e$ and
$$
a_{n+1}=\frac1e+\frac1e\sum_{k=1}^\infty\frac{a_n^k}{k!}=e^{a_n-1}\;\;(n=1,2,\dots).$$
Evidently, the sequence $\{a_n\}$ increases and $a_n\in[0,1].$ Since the function $f(x):=\,e^{x-1}-x$ decreases on $[0,1],$ the function
$e^{x-1}$ has the only fixed point $x=1$. Hence, $\lim_{n\to\infty}a_n=1,$ which proves the lemma.
\end{proof}

\begin{proof}[Proof of Theorem \ref{psi theorem}]
Consider the functions $h_n=(K^n1)^*,$ $n\geq0.$ Since the operator $K$ maps equimeasurable functions to equimeasurable ones, we have
\begin{equation}\label{recurrence}
(Kh_n)^*=h_{n+1}.
\end{equation}
By Lemma \ref{power of operator K}, $\mes({\rm supp}\,h_n)\to0$ as $n\to\infty.$ Hence, the series
\begin{equation}\label{psi definition}
g_{\varepsilon}=\sum_{n=0}^{\infty}\varepsilon^nh_n
\end{equation}
converges everywhere on the interval $(0,1]$ for every $\varepsilon>0$ and the function $g_{\varepsilon}$ decreases. Moreover, it follows from the definition of the operator $K$ (see \eqref{40-(-1)}) that $\|K\|_{L_1}=1.$ Hence, if $0<\varepsilon<1,$ then the series \eqref{psi definition} converges in $L_1$ and $g_\varepsilon\in L_1.$ We shall show that the assertions of the theorem hold for the family $\{M_{\psi_{\varepsilon}}\}_{\varepsilon>0}$, where $\psi_{\varepsilon}(t)=\int_0^t g_{\varepsilon}(s)\,ds$ $(0\le t\le1).$

1. Let us prove that the operator $K$ is bounded in $M_{\psi_{\varepsilon}}.$ The extreme points of the unit ball in this space are equimeasurable with  $g_{\varepsilon}$ \cite{Ryff} and, therefore, it is sufficient to show that $Kg_{\varepsilon}\in M_{\psi_{\varepsilon}}.$ Since $K$ is bounded in $L_1,$ then
$$Kg_{\varepsilon}=\sum_{n=0}^{\infty}\varepsilon^nKh_n\prec\sum_{n=0}^{\infty}\varepsilon^nh_{n+1}\leq\frac1{\varepsilon}\sum_{n=0}^{\infty}\varepsilon^nh_n=\frac1{\varepsilon}g_{\varepsilon}.$$
Here, the first inequality follows from \eqref{recurrence} and the well-known property of Hardy-Littlewood submajorization (see, for example, \cite[\S\,2.2]{lit-9}). Thus, $Kg_{\varepsilon}\in M_{\psi_{\varepsilon}}.$

2. Now assume that $E\in\mathbb{K}$. As we mentioned earlier, this assumption guarantees that $C=||K||_{E\to E}<\infty$. Evidently, $||h_n||_E\leq C^n||1||_E.$ Therefore, for every $\varepsilon<C^{-1}$ the series \eqref{psi definition} converges in $E$ and $g_{\varepsilon}\in E.$ Since the space $E$ is either separable or $E=E^{\prime\prime}$, we have that $x\in E$ and $y\prec x$ imply that $y\in E$ and $||y||_E\leq||x||_E.$ Hence, the unit ball of the space $M_{\psi_{\varepsilon}}$ is a subset of $E.$ Therefore, $M_{\psi_{\varepsilon}}\subset E.$

3. Let the function $g_{\varepsilon}$ be as in \eqref{psi definition} and let $0<\varepsilon<\delta.$ Arguing as in the proof of the Theorem  7.2 in \cite{lit-3}, one can obtain
$$\lim_{t\to 0}\frac{h_{n+1}(t)}{h_{n}(t)}=\infty.$$
Therefore, for every $m=1,2,\dots$
\begin{eqnarray*}
\limsup_{t\to 0}\frac{g_{\varepsilon}(t)}{g_{\delta}(t)} &=&
\limsup_{t\to
0}\left(\sum_{n=1}^{\infty}\varepsilon^nh_n(t)\right)\cdot\left(
\sum_{n=1}^{\infty}\delta^nh_n(t)\right)^{-1}=\dots\\
\dots & = &\limsup_{t\to
0}\left(\sum_{n=m}^{\infty}\varepsilon^nh_n(t)\right)\cdot\left(
\sum_{n=m}^{\infty}\delta^nh_n(t)\right)^{-1}\le
\left(\frac\varepsilon\delta\right)^m.
\end{eqnarray*}
Therefore, $\lim_{t\to 0}\frac{g_{\varepsilon}(t)}{g_{\delta}(t)}=0$ and the assertion \eqref{non-equivalent} follows immediately.

4. According to the introduction, the operator $K$ acts boundedly in the space $(\exp L_1)_0.$ Hence, the fourth assertion follows from the second and third ones.
\end{proof}

Let $\varphi_n(t):=\int_0^th_n(s)\,ds$ $(0\le t\le 1)$ and let
$M_{\varphi_n}$ be the corresponding Marcinkiewicz space. We have, $M_{\varphi_n}\subset M_{\varphi_{n+1}}\subset(\exp L_1)_0$ $(n=1,2,\dots)$ and so in a certain sense the spaces $M_{\varphi_n}$, $n\ge 1$ may be viewed as ``approximations'' of the space $(\exp L_1)_0$. By \cite[Theorem 7.2]{lit-3}, we have $M_{\varphi_n}\subset E$ for every r.i. space $E\in \mathbb{K}$ and every $n=1,2,\dots$ This suggests a rather natural conjecture that $(\exp L_1)_0$ is the minimal r.i. space with the Kruglov property. However, the following consequence from Theorem \ref{psi theorem} shows that the class of r.i. spaces with the Kruglov property has no minimal element.

\begin{cor} \label{corol: there isn't minim. r.i.s.} For every r.i. space $E\in\mathbb{K}$ there exists an r.i. space $F\in\mathbb{K}$ such that $F\subsetneqq E.$
\end{cor}

Contrary to the case of Marcinkiewicz spaces, all Lorentz spaces with the Kruglov property lie ``on the one side'' of the space $\exp L_1.$

\begin{thm}\label{Lorentz theorem} Let $\varphi$ be an increasing concave function on the interval $[0,1]$ such that $\varphi(0)=0.$ If $\Lambda_\varphi\in\mathbb{K},$ then $\Lambda_\varphi\supset \exp L_1.$
\end{thm}

Let us prove the following Lemma first.
\begin{lem}\label{embedding of Lor. sp.} Let $\varphi$ be an increasing function on the interval $[0,1]$ and let $\varphi(0)=0.$ If $\varphi$ satisfies condition \eqref{40-2}, then
\begin{equation}\label{equa1}
 \sum_{k=1}^\infty \varphi(2^{-k})\le A\varphi(1).
\end{equation}
Here, $A>0$ depends only on $M$ from \eqref{40-2}.
\end{lem}
\begin{proof}
According to \eqref{40-2}, for every $i\in\N$
$$
\sum_{j=1}^\infty \varphi(2^{-ij}j^{-j})\le M\varphi(2^{-i})$$
or, equivalently,
\begin{equation}\label{equa3}
\sum_{j=1}^\infty \varphi(2^{-j(i+[\log_2j])})\le M\varphi(2^{-i}).
\end{equation}
Straightforward calculations show that the quantity
$$
\alpha_n:=|\{(i,j)\in \N^2:\,j(i+[\log_2j])\le n\}|$$
satisfies the condition $\lim_{n\to\infty}n^{-1}\alpha_n=\infty.$ Hence,  $\alpha_n\geq
(M+1)n$ for some $m\in\N$ and for every $n\geq m.$ It follows from \eqref{equa3} and the monotonicity of $\varphi$ that for every $l>m$
$$
(M+1)\sum_{n=m}^l\varphi(2^{-n})\le \sum_{i=1}^l \sum_{j=1}^\infty
\varphi(2^{-j(i+[\log_2j])})\le M\sum_{i=1}^l\varphi(2^{-i}).
$$
Thus,
$$
\sum_{n=m}^l\varphi(2^{-n})\le M\sum_{i=1}^{m-1}\varphi(2^{-i}).
$$
Note that $m$ depends only on $M$ and not on $\varphi,$ while $l>m$ is arbitrary. The inequality \eqref{equa1} follows immediately.
\end{proof}

\begin{proof}[Proof of Theorem \ref{Lorentz theorem}]
According to the introduction, condition \eqref{40-2} is equivalent to the condition $\Lambda_\varphi\in\mathbb{K}$ \cite{lit-3}. Therefore,
Lemma \ref{embedding of Lor. sp.} implies that condition \eqref{equa1} holds. Moreover, by \cite{lor}, we have
$$
\|x\|_{\exp L_1}\asymp \sup_{0<t\le 1}x^*(t)\log_2^{-1}(2/t)$$
 and therefore to prove the embedding $\Lambda_\varphi \supset\exp L_1$ it is sufficient to prove only that $\log_2(2/t)\in \Lambda_\varphi.$ The latter follows from the following estimates:
\begin{eqnarray*}
\|\log_2(2/t)\|_{\Lambda_\varphi}
&=&\int_0^1\log_2(2/t)\,d\varphi(t)=\sum_{k=1}^\infty
\int_{2^{-k}}^{2^{-k+1}}\log_2(2/t)\,d\varphi(t)\\
&\le&
\sum_{k=1}^\infty(k+1)(\varphi(2^{-k+1})-\varphi(2^{-k}))=2\varphi(1)+
\sum_{k=1}^\infty\varphi(2^{-k})<\infty.
\end{eqnarray*}
\end{proof}
\skip 0.2cm

{\centering\section{Estimates of distribution functions}}

We will use the following approximation of $Kf,$ where $f$ is an arbitrary measurable function on the interval $[0,1].$

Let $m\in\N,$ $g_m(t)=\sigma_{\frac1m}f(t)$ and let $\{h_{m,i}\}_{i=1}^m$ be independent functions equimeasurable with $g_m.$ The sequence
\begin{equation}\label{40-3.1}
    H_mf(t)=\sum_{i=1}^m h_{m,i}(t)\;\;(0\le t\le 1)
\end{equation}
weakly converges to $Kf$ when $m\to\infty$ (in the sense of convergence of distribution functions) (see \cite[1.6, p.~11]{lit-2}) or \cite[Theorem~3.5]{lit-3}).

In particular, if $n\in\N,$ $a_k\geq0$ $(1\leq k\leq n)$ and
\begin{equation}\label{40-3}
    f_a(t)=\sum_{k=1}^n a_k\chi_{\left(\frac{k-1}n, \frac kn\right)}
    (t)\;\;(0\le t\le 1),
\end{equation}
then
\[
 g_m(t)=\sigma_{\frac1m}f_a(t)=\sum_{k=1}^n a_k\chi_{\left(\frac{k-1}{nm},
 \frac{k}{nm}\right)}(t)\;\;(m\in\N).
\]
In the latter case, we denote
\begin{equation}\label{40-4}
    H_ma(t):=H_mf_a(t)=\sum_{i=1}^m h_{m,i}(t).
\end{equation}
In addition, let $Ch(r)$ be the number of permutations $\pi$ of the set $\{1,2,\ldots,r\}$ such that $\pi(i)\not=i$ for every $i=1,2,\ldots,r.$ It is well known (see \cite[p.~20]{lit-10}) that
\begin{equation}\label{40-5}
 \frac13r!\leq Ch(r)\leq r!\;\;(r\in\N).
\end{equation}

We are going to compare distribution functions of $H_ma$ and $T_{nm}b$, where
\[
 b=(\underbrace{a_1,a_1,\ldots,a_1}_m,\underbrace{a_2,a_2,\ldots,a_2}_m,\ldots,
 \underbrace{a_n,a_n,\ldots,a_n}_m).
\]
\begin{lem}\label{lem-40-1}
For every $n, m\in\N$ and every $\tau>0$
\[
 \mes\{t:\, H_ma(t)>\tau\}\leq
 3\mes\{t:\,T_{nm}b(t)>\tau\}.
\]
\end{lem}
\begin{proof} The function $H_ma(t)$ (respectively, $T_{nm}b(t)$) only takes values of the form $\displaystyle\sum_{i=1}^nk_ia_i$, where $k_i\in\mathbb{Z},$ $k_i\geq0$ for all $i=1,2,\ldots,n$ and $\displaystyle\sum_{i=1}^n k_i\leq m$ (respectively, $\displaystyle\sum_{i=1}^n k_i\leq mn$). Therefore, it is sufficient to prove that
\[
 \mes\left\{
  t:\,H_ma(t)=\sum_{i=1}^n k_ia_i
 \right\}\leq3
 \mes\left\{
  t:\,T_{nm}b(t)=\sum_{i=1}^n k_ia_i
 \right\}
\]
for any choice of $k_i\in\N$, $\displaystyle \sum_{i=1}^n k_i=q\leq m$. Note, that it is sufficient to consider only the case when
$$
\sum_{i=1}^n k_ia_i\ne \sum_{i=1}^n k_i'a_i\;\;\mbox{provided that}\;\;(k_1,k_2,\dots,k_n)\ne (k_1',k_2',\dots,k_n').
$$
Hence, $H_ma(t)$ equals $\sum_{i=1}^n k_ia_i$ if and only if exactly $k_i$ (respectively, $m-q$) of the functions $h_{m,j}(t)$ $(j=1,\dots,m)$ take the value $a_i$ (respectively, 0). Since the functions $h_{m,j}$ are independent, we obtain
\begin{eqnarray}
\mes\left\{t:\ H_ma(t)=\sum_{i=1}^nk_ia_i
\right\}&=&C_m^{m-q,k_1,\cdots,k_n} \left(1-\frac1m
\right)^{m-q}\left(\frac1{mn} \right)^{k_1+\cdots+k_n}\nonumber\\
&\le &
C_m^{m-q,k_1,\cdots,k_n}\left(\frac1{mn}\right)^q,\label{39-5}
\end{eqnarray}
where
$$
 C_m^{m-q,k_1,\cdots,k_n}=
\frac{m!}{(m-q)!k_1!\ldots k_n!}.$$

On the other hand, it follows from \eqref{40-1} and \eqref{40-5} that
\begin{eqnarray*}
 \mes\left\{t:\, T_{mn}b(t) = \sum_{i=1}^n k_ia_i
 \right\}&=& C_m^{k_1}C_m^{k_2}\ldots C_m^{k_n}Ch(mn-q)
 \frac1{(mn)!}\\
 &\ge&
\frac{(m!)^n(mn-q)!}{3(m-k_1)!\cdots(m-k_n)!k_1!\cdots k_n!(mn)!}.
\end{eqnarray*}
Since
$$(m-k_1)!\cdots (m-k_n)!\leq (m!)^{n-1}(m-q)!$$
and
$$\frac{(mn-q)!}{(mn)!}\geq\frac1{(mn)^q},$$
we have
\begin{eqnarray*}
 \mes\left\{t:\, T_{mn}b(t)
= \sum_{i=1}^n k_ia_i \right\}&\ge &
\frac{m!(mn-q)!}{3k_1!\cdots k_n!(m-q)!(mn)!}\\
&\ge & \frac{m!}{3(m-q)!k_1!\cdots k_n!}\cdot\frac1{(mn)^q}.
\end{eqnarray*}
The assertion follows now from this inequality and inequality \eqref{39-5}.\end{proof}

\begin{lem}\label{lemma-40-2} If $n,k\in\N$, $n\geq 4,$ $k\le n,$ then
\begin{equation*}\label{40-7}
    \frac{(n-k)!}{n!}\leq 2\frac{(k-1)!}{n^k}.
\end{equation*}
\end{lem}
\begin{proof} Since $j(n-j)>n$ for $2\le j\le n-2,$ we have
$$\frac{n^k(n-k)!}{n!(k-1)!}=\prod_{j=1}^{k-1}\frac{n}{j(n-j)}\leq\left(\frac{n}{n-1}\right)^2<2.$$
\end{proof}
\skip 0.2cm

 Now we continue the study begun in Lemma ~\ref{lem-40-1} of the connections between the distribution functions of $T_na$ and $H_ma.$  Whereas the estimate obtained in Lemma~\ref{lem-40-1} holds for every $m$ and $n$, the converse inequality holds only asymptotically when $m\to\infty.$

\begin{lem}\label{lemm-40-3} Let $n\in\N$, $a=(a_1,a_2,\ldots,a_n)\geq0$, $\tau>0$. For every sufficiently large $m\in N,$ the following inequality is valid:
\[
 \mes\{t:\, T_na(t)>\tau\}\leq 12 \mes\{t:\,2H_ma(t)>\tau\}.
\]
\end{lem}
\begin{proof}
Assume first that $n\geq 4.$ Let $A=\{1,2,.,n\}.$ Denote
$$S(U):=\sum_{j\in U}a_j$$
for every $U\subset A.$ Without loss of generality, we may assume that $n=2s$ $(s\in\N)$, $a_i>0$ and $S(U_1)\ne S(U_2)$ if $U_1\ne U_2.$ Denote by $\mathcal{A}_i$ the collection of all sets $U\subset A$ with $|U|=i$ $(i=1,2,\cdots,n).$ Hence, $\mathcal{A}=\cup_{i=1}^n\mathcal{A}_i$ is the collection of all non-empty subsets of the set $A.$ Let us represent the set $\mathcal{A}$ in another way.

Let $U\in\mathcal{A}_k$ for some $k=1,2,\cdots,s.$ Denote $\mathcal{A}_U$ (respectively, $\mathcal{B}_U$) the collection of all sets $V\subset A$ such that $V\supset U,$ $V\in\mathcal{A}_{2k}$ (respectively, $V\in\mathcal{A}_{2k-1})$ and $S(V\setminus U)\le
S(U).$ Since
$$
\bigcup_{U\in\mathcal{A}_k}\mathcal{A}_U=\mathcal{A}_{2k}\;\;\mbox{and}\;\;
\bigcup_{U\in\mathcal{A}_k}\mathcal{B}_U=\mathcal{A}_{2k-1}\;\;(k=1,2,.,s),$$
then
\begin{equation}\label{40-7.1}
\mathcal{A}=\bigcup_{k=1}^s\bigcup_{U\in\mathcal{A}_k}\big(\mathcal{A}_U\cup\mathcal{B}_U\big).
\end{equation}
It follows from the definition of $\mathcal{A}_U$ and $\mathcal{B}_U$ that for every $V\in\mathcal{A}_U\cup\mathcal{B}_U$
\begin{equation}
\label{40-7.2} S(U)\le S(V)\le 2S(U).
\end{equation}

Note that $T_na(t)$ is a step function with values of the form $S(V),$ where $V\in\mathcal{A}.$ If $|V|=r,$ then \eqref{40-5} implies that
$$
\mes\{t:\,T_na(t)=S(V)\}=\frac{Ch(n-r)}{n!}\leq\frac{(n-r)!}{n!}.$$
Also, if $|U|=k$ $(k=1,2,.,s),$ then
$$|\mathcal{A}_U|\le C_{n-k}^k=\frac{(n-k)!}{k!(n-2k)!}$$
and similarly
$$|\mathcal{B}_U|\le C_{n-k}^{k-1}=\frac{(n-k)!}{(k-1)!(n-2k+1)!}.$$
Therefore, \eqref{40-7.1} and \eqref{40-7.2} imply that
\begin{eqnarray}
\mes\{t:\,T_na(t)>\tau\} &\le &\sum_{k=1}^s\sum_{U\in \mathcal{A}_k}\Bigg(
\sum_{V\in\mathcal{A}_U,S(V)>\tau}\mes\{t:\,T_na(t)=S(V)\}\nonumber\\
&+&
\sum_{V\in\mathcal{B}_U,S(V)>\tau}\mes\{t:\,T_na(t)=S(V)\}\Bigg)\nonumber\\
&\le &\sum_{k=1}^s\sum_{U\in\mathcal{A}_k,S(U)>\tau/2}\Bigg(\frac{(n-2k)!}{n!}\cdot\frac{(n-k)!}{k!(n-2k)!}\nonumber\\
&+&
\frac{(n-2k+1)!}{n!}\cdot\frac{(n-k)!}{(k-1)!(n-2k+1)!}\Bigg)\nonumber\\
&\le& 2\sum_{k=1}^s\sum_{U\in\mathcal{A}_k,S(U)>\tau/2}\frac{(n-k)!}{(k-1)!n!}. \label{40-7.3}
\end{eqnarray}

Let us now estimate the distribution function of $H_ma(t)$ from below. For every $U\in\mathcal{A}_k,$ $S(U)>\tau/2,$ let $F_U$ be the set of all $t\in [0,1]$ such that there exists a set $W\subset\{1,2,\cdots,m\}$ and a bijection $\sigma:\,W\to U,$ such that $|W|=k$ (we assume that $m\ge n$) and
$h_{m,j}(t)=a_{\sigma(j)}$ if $j\in W,$ and $h_{m,j}(t)=0$ if $j\not\in W.$ Thus, for $t\in F_U$
\begin{equation}\label{40-8}
H_ma(t)=\sum_{j=1}^m h_{m,j}(t)=S(U)>\frac{\tau}{2}.
\end{equation}
The independence of the functions $h_{m,j}(t)$ $(j=1,2,\cdots,m)$ implies
\begin{eqnarray*}
\mes (F_U) &=& C_m^{k}k! \frac1{(mn)^{k}}
    \Big(1-\frac1{m}\Big)^{m-k}\\
   & = & \frac{m(m-1)\cdot\dots\cdot(m-k+1)}{m^k}\cdot
    \Big(1-\frac1{m}\Big)^{m-k}\cdot\frac1{n^k}.
\end{eqnarray*}
Since
\[
 \lim_{m\to\infty}\frac{m(m-1)\ldots(m-k+1)}{m^{k}}=1
\]
and
\[
 \lim_{m\to\infty}\Big(1-\frac1{m}\Big)^{m-k}
 =\frac1e >\frac13,
\]
we obtain
\begin{equation}\label{40-8.1}
 \mes (F_U) > \frac13\cdot\frac1{n^k}
\end{equation}
for all sufficiently large $m\in\N$ and for all $k\le s.$

Note that $F_U\cap F_{U'}=\emptyset$ if $U\neq U'.$ Indeed, let $i\in U\setminus U'.$ For every $t\in F_U$ there exists $j\in\{1,2,.,m\}$ such that $h_{m,j}(t)=a_i.$ However, if $t\in F_{U'},$ then either $h_{m,j}(t)=a_l\ne a_i$ or $h_{m,j}(t)=0\ne a_i.$ Hence, equations \eqref{40-8.1} and \eqref{40-7.3} and Lemma \ref{lemma-40-2} imply that
\begin{eqnarray*}
\mes\{t:\,2H_ma(t)>\tau\} &=& \sum_{k=1}^s\sum_{U\in\mathcal{A}_k,S(U)>\tau/2}\mes(F_U)\\ &\ge & \frac13 \sum_{k=1}^s\sum_{U\in\mathcal{A}_k,S(U)>\tau/2}\frac1{n^k}\\ &\ge & \frac16
\sum_{k=1}^s\sum_{U\in\mathcal{A}_k,S(U)>\tau/2}\frac{(n-k)!}{(k-1)!n!}\\&\ge & \frac1{12}
\mes\{t:\,T_na(t)>\tau\}.
\end{eqnarray*}
This estimate proves the lemma for $n\ge 4.$

If $1\leq n<4,$ then it is easy to show (see the argument preceding equation \eqref{40-8.1}) that
\begin{equation*}
\mes\{t:\, T_na(t)>\tau\}\leq 5\mes\{t:\,2H_ma(t)>\tau\}
\end{equation*}
for all sufficiently large $m\in\N$ and every $\tau>0$.
\end{proof}

\begin{rem} The estimate
\[
 \mes\{t:\,T_na(t)>\tau\}\leq C\mes\{t:\,H_na(t)>\tau\}\;\;(\tau>0)
\]
fails for any constant $C$ independent of $n\in\N$. Indeed, if $a_1=a_2=\ldots=a_n=1,$ then
\[
 \mes\{t:\,T_na(t)=n\}=\frac1{n!},
\]
while
\[
 \mes\{t:\, H_na(t)=n\}=\frac1{n^n}.
\]
\end{rem}
\vskip 0.2cm

{\centering\section{The Kruglov property and random permutations}}

\begin{thm}\label{teo-40-6} Let $E$ be an r.i. space. The operator $K$ acts boundedly on $E$ if and only if the sequence of operators $T_n$ is uniformly bounded in $E.$
\end{thm}
\begin{proof}
We are going to use notations \eqref{40-1}, \eqref{40-3} and \eqref{40-4}.

Necessity. It follows from Lemma \ref{lemm-40-3} that for arbitrary $n\in\N$, $a=(a_1,a_2,\ldots,a_n)\geq0$, $\tau>0$ and every sufficiently large
$m\in\N$ we have
\[
 \mes\{t:\, T_na(t)>\tau\}\leq 12 \mes\{t:\,2H_ma(t)>\tau\}.
\]
As we pointed out in the preceding section, $H_ma\Rightarrow Kf_a$ when $m\to\infty.$ Therefore, \cite[\S\,6.2]{lit-11},
$$
\mes\{t:\,H_ma(t)>\tau\}\to\mes\{t:\,Kf_a(t)>\tau\}\;\;(m\to\infty)
$$
if the right-hand side is continuous at $\tau>0.$ Hence, the convergence is valid for all but countably many values of $\tau.$ Hence, for all such $\tau,$ we have
\[
 \mes\{t:\, T_na(t)>\tau\}\leq 12 \mes\{t:\,2Kf_a(t)>\tau\}.
\]
Both functions in the last inequality are monotone and right-continuous. Therefore, this inequality holds for every $\tau>0.$

It is well known (see \cite[\S\,2.4.3]{lit-9}), that for every r.i. space $E$ the relation $y\in E$ together with the inequality
$$\mes\{t:\,|x(t)|>\tau\}\leq C\mes\{t:\,|y(t)|>\tau\}\;\;(\tau>0)$$
imply that $x\in E$ and $\|x\|_E\leq \max(C,1)\|y\|_E$. Therefore, by the preceding inequality
$$\|T_nf_a\|_E\le 24\cdot\|Kf_a\|_E$$
or
$$
\sup\{\|T_nf_a\|_E:\,\|f_a\|\le 1\}\le 24\cdot\|K\|_E.$$
By the definition of the operator $T_n,$ we have $T_nx=T_nf_{a_n(x)},$ where $a_n(x)=(a_{n,k}(x))_{k=1}^n,$ $a_{n,k}(x)=n\int_{\frac{k-1}{n}}^{\frac{k}{n}}x(s)\,ds.$ Since $\|f_{a_n(x)}\|_E\le \|x\|_E$ \cite[\S\,2.3.2]{lit-9} and due to the assumption that $E$ is either separable or coincides with its second K\" othe dual, we obtain
\[
 \sup_n\|T_n\|_E\leq 24\cdot\|K\|_E.
\]

Sufficiency. Assume that $\displaystyle\sup_n\|T_n\|_E=C<\infty$. It follows from Lemma~\ref{lem-40-1} and \cite[\S\,2.4.3]{lit-9} that
\[
 \|H_mf_a\|_E\leq3\|T_{nm}\|_E\|f_a\|_E\leq3C\|f_a\|_E.
\]
Since $H_mf_a\Rightarrow Kf_a$ when $m\to\infty,$ it follows from \cite[Proposition~1.5]{lit-2} that
\begin{equation}\label{40-8.2}
 \|Kf_a\|_{E''}\leq 3C\|f_a\|_E.
\end{equation}

Let now $f=f^*\in E$ be arbitrary. If
$$f_n(t)=\sum_{k=1}^{2^n} f(k2^{-n})\chi_{((k-1)2^{-n},k2^{-n})}(t)\;\;(0\le t\le 1),\;\;n\in\N,$$
then $f_n(t)\uparrow f(t)$ a.e., and, therefore, $f_n\Rightarrow f$ \cite[\S\,6.2]{lit-11}. If $\varphi_n$ and $\varphi$ are the characteristic functions of $f_n$ and $f$ respectively, then $\varphi_n(t)\to \varphi(t)$ $(t\in\mathbb{R})$ (\cite[\S\,6.4]{lit-11}). In view of \cite[1.6]{lit-2}, we have
$$\varphi_{K\xi}(t)=\exp(\varphi_\xi(t)-1)$$
for every random variable $\xi.$ Hence, $\varphi_{Kf_n}(t)\to\varphi_{Kf}(t)$ $(t\in\mathbb{R}),$ i.e.  $Kf_n\Rightarrow Kf$. Thanks to \eqref{40-8.2}, we have
$$\|Kf_n\|_{E''}\leq 3C\|f_n\|_E\le 3C\|f\|_E\;\;(n\in\N).$$
Thus, using \cite[Proposition~1.5]{lit-2} once more, we obtain
$$\|Kf\|_{E''}\leq 3C\|f\|_E.$$
Since the distribution function of $Kf$ depends only on the distribution function of $f$, it follows from the preceding inequality that
the operator $K$ boundedly maps $E$ into $E''$. If $E=E'',$ then we are done. It remains to consider the case when $E\neq E''$. In this case, the space $E$ is separable. First of all, using the fact that every function $f\in E''$, $f\ge 0,$ is the a.e. limit of its truncations
$\tilde{f}_n:=f\chi_{\{f_n\le n\}}$ $(n\in\N)$ and arguing as above, one can infer that the operator $K$ acts boundedly in $E''.$ Therefore, by \cite[Theorem~7.2]{lit-3}, the function
$$g(t):=\frac{\ln(e/t)}{\ln(\ln(\ln(a/t)))},$$
where $a>0$ is sufficiently large, belongs to $E''.$ Now, if
$$\psi(u):=\frac{u\ln(e/u)}{\ln(\ln(\ln(a/u)))}\;\;(0<u\le 1),$$
then the Marcinkiewicz space $M_\psi\subset E''$. Hence, in view of separability of the space $E$, we have
$$(M_\psi)_0\subset (E'')_0=E_0=E.$$
It is easy to check that
$$h(t):=\frac{\ln(e/t)}{\ln(\ln(a/t))}\in(M_\psi)_0,$$ whence, $h\in E.$ This and \cite[Th.~4.4]{lit-3} imply that
\begin{equation}\label{40-8.3}
K:\,L_\infty\to E.
\end{equation}

Let now $f\in E.$ Since $E$ is separable, there exists a sequence $\{f_n\}\subset L_{\infty}$ such that $||f_n-f||_E\to 0.$ Since $K:E\to E'',$ we have $\|Kf_n-Kf\|_{E''}\to 0.$ On the other hand, by \eqref{40-8.3} and taking into account that the embedding $E\subset E''$ is isometric, we have $\{Kf_n\}\subset E,$  whence $Kf\in E.$
\end{proof}

\begin{rem} It follows from the proof above that the following estimate holds in every r.i. space $E$
\[
 \frac1{24}\sup_n\|T_n\|_E\leq\|K\|_E\leq3\sup_n\|T_n\|_E.\qed
\]
\end{rem}
\vskip 0.2cm

We are going to infer some corollaries from Theorem~\ref{teo-40-6}. Let $n\in\N$ and let $S_n$ be the set of all permutations of the set $\{1,2,\ldots,n\}.$ Fix a map $l=l_n$ from $S_n$ onto the set $\{1,2,\ldots,n!\}$. Recall that the earlier definition of the operator $A_n$ acting from $\R$ is given by \eqref{40-1}. We are now in a position to extend this definition to the set of matrices $x=(x_{i,j})_{1\le i,j\le n}$ as follows
\begin{equation*}
    A_nx(t)=\sum_{i=1}^n x_{i,\pi(i)},\quad t\in\left(
    \frac{l(\pi)-1}{n!},\frac{l(\pi)}{n!}
    \right).
\end{equation*}
One of the major results of \cite{lit-7} (see Corollary~8 there) says that if the sequence of operators $\{A_n\}_{n\ge 1}$ is uniformly bounded on the set of diagonal matrices, then it is uniformly bounded on the set of all matrices. Applying Theorem~\ref{teo-40-6}, we obtain
\begin{cor}\label{corol-40-9} If an r.i. space $E\in\mathbb{K},$ then for every $n\in\N$ and every $x=(x_{i,j})_{1\le i,j\le n}$
$$
\|A_nx\|_E\le C\left(\Big\|\sum_{k=1}^nx_k^*\chi_{\left(\frac{k-1}n,
\frac{k}{n}\right)}\Big\|_E+\frac1n\sum_{k=n+1}^{n^2}x_k^*\right).$$
Here, $(x_k^*)_{k=1}^{n^2}$ is a decreasing permutation of the sequence $(|x_{i,j}|)_{i,j=1}^n$ and $C>0$ does not depend either on $n$ or $x.$
\end{cor}

\begin{cor}\label{corol-40-8} The operators $T_n$, $n\ge 1$ are uniformly bounded in Orlicz space $\exp L_p$ if and only if $p\leq1.$
\end{cor}

Indeed, the Orlicz space $\exp L_p$ has the Kruglov property if and only if $p\leq 1$ (see \cite[2.4, p.~42]{lit-2}). The preceding corollary now follows immediately from Theorem~\ref{teo-40-6}.

Theorem~\ref{teo-40-6} and Corollary~\ref{corol: there isn't minim. r.i.s.} imply
\begin{cor} If $E$ is an r.i. space and if $\sup_n\|T_n\|_E<\infty,$ then there exists an r.i. space $F\subsetneq E,$ such that $\sup_n\|T_n\|_F<\infty.$
\end{cor}

If $E$ is an r.i. space and $p\geq1,$ then $E(p)$ denotes the space of all measurable functions $x$ on the interval $[0,1]$ such that $|x|^p\in E.$ We equip $E(p)$ with the norm
$$\|x\|_{E(p)}=\|\,|x|^p\,\|_E^{1/p}.$$
It is well known that $E(p)\subset E$ and $\|x\|_E\le \|x\|_{E(p)}$ for all $x\in E(p)$  \cite[1.d]{lit-8}.

Let $E$ and $F$ be r.i. spaces such that $E\subset F$ and $K:E\to E.$ This does not imply in general that $K:F\to F$ \cite[Corollaries 5.6 and 5.7]{lit-3}. However, we have
\begin{cor} If the operator $K$ is bounded in $E(p),$ then it is bounded in $E.$
\end{cor}
\begin{proof} By Theorem \ref{teo-40-6}, it is sufficient to prove that the uniform boundedness of operators $T_n$, $n\ge 1$ in $E(p)$ implies the uniform boundedness of operators $T_n$, $n\ge 1$ in $E$.

Let $x=(x_1,x_2,\dots,x_n)\in \mathbb{R}^n,$ $x\ge 0$ and $\|T_nx\|_{E(p)}\le C\|x\|_{E(p)}$ $(n\in\N).$ It means that,
$$
\|(T_nx)^p\|_E^{1/p}\le C\|x^p\|_E^{1/p}.$$
If $x^p=y,$ then
$$
\|(T_ny^{1/p})^p\|_E\le C^p\|y\|_E.$$
It follows from the definition of the operator $T_n$, $n\ge 1$ that $(T_ny^{1/p})^p\ge T_ny,$ Hence, $\|T_ny\|_E\le C^p\|y\|_E$, $n\ge 1$. Thus, the operators $T_n$, $n\ge 1$ are uniformly bounded in $E$.
\end{proof}

\begin{center}
Astashkin S.V.\\[2mm]
Samara State University\\[2mm]
443011 Samara, Acad. Pavlov, 1 \\[2mm]
Russian Federation\\[2mm]
e-mail: {\it astashkn@ssu.samara.ru}
\vskip 0.4cm

Zanin D.V.\\[2mm]
School of Computer Science, Engineering and Mathematics\\[2mm]
Flinders University, Bedford Park, SA 5042 Australia\\[2mm]
e-mail: {\it zani0005@infoeng.flinders.edu.au}
\vskip 0.4cm

Semenov E.M.\\[2mm]
Voronezh State University\\[2mm]
394006, Voronezh, University pl., 1\\[2mm]
Russian Federation\\[2mm]
e-mail: {\it semenov@func.vsu.ru}
\vskip 0.4cm

Sukochev F.A.\\[2mm]
School of Mathematics and Statistics\\[2mm]
University of New South Wales\\[2mm]
Kensington NSW 2052 Australia\\[2mm]
e-mail: {\it f.sukochev@unsw.edu.au}
\end{center}

%
%
%
%
%
%
\end{document}